\documentclass[12pt]{amsart}
\usepackage{amssymb, longtable}
\usepackage[colorlinks=true,citecolor=black,linkcolor=black,urlcolor=blue]{hyperref}

\theoremstyle{plain}
\newtheorem{theorem}{Theorem}
\newtheorem{proposition}[theorem]{Proposition}
\theoremstyle{definition}
\newtheorem{definition}[theorem]{Definition}

\newcommand{\ZZ}{\mathbb{Z}}

\newenvironment{numberlists}[1][3\parindent]
 {\begin{list}{}{%
   \leftmargin=#1\relax
   \rightmargin=\leftmargin
   \itemsep=\jot
   \parsep=0pt
   \partopsep=0pt
   \labelsep=0pt}}
 {\end{list}}
\newcommand\numlist[2]{%
  \item[]\makebox[0pt][r]{$\{$}%
  \begingroup
  \begingroup\lccode`~=`,\lowercase{\endgroup\def~}{\mathcomma\penalty0 }%
  \mathcode`,="8000
  \thinmuskip=6mu plus 6mu minus 2mu
  $#1\}.$%
  \endgroup
}
\mathchardef\mathcomma=\mathcode`,

\title{All $(96,20,4)$ difference sets and related structures}
\author{Omar A. AbuGhneim}
\address{Department of Mathematics, University of Jordan, Amman, Jordan}
\email{o.abughneim@ju.edu.jo}
\author{Dylan Peifer}
\address{The Department of Mathematics, Cornell University, Ithaca, New York 14853}
\email{djp282@cornell.edu}
\author{Ken W. Smith}
\address{The Department of Mathematics and Statistics, Sam Houston State University, Huntsville, Texas 77340}
\email{kenwsmith@shsu.edu}

\begin{document}

\maketitle

\begin{abstract}
In 1978, Robert Kibler \nocite{Kibler1978} at the National Security Agency in Fort Meade, Maryland published a description of all noncyclic difference sets with $k < 20$.
Kibler's decision to stop his extensive computer search for difference sets at block size 19 was motivated partly by the difficult barrier at $k=20$, the difference sets with parameters $(96,20,4)$.
In this paper, we announce the completion of the search for all $(96,20,4)$ difference sets, relying on the computer software {\tt GAP} and the work of numerous authors over the last few decades.
The difference sets and the symmetric designs they create are summarized and links are provided to webpages which explicitly list the difference sets.
In addition, we use these $(96,20,4)$ difference sets to construct all $(96, 20, 4, 4)$ and $(96, 19, 2, 4)$ partial difference sets and briefly look at the corresponding strongly regular graphs.
\end{abstract}

\section{Introduction}\label{sec:intro}

A $(v, k, \lambda)$ difference set is a subset $D$ of size $k$ in a group $G$ of order $v$ with the property that for every nonidentity $g \in G$, there are exactly $\lambda$ ordered pairs $(x, y) \in D \times D$ such that $xy^{-1} = g$.
One may identify the set $D$ with an element $\hat{D}$ in the group ring $\ZZ[G]$.
In this case write
\[
\hat{D} = \sum_{g \in D}g, \qquad \hat{D}^{(-1)} = \sum_{g\in D}g^{-1}, \qquad \hat{G} = \sum_{g\in G} g
\]
and then $D$ is a difference set if the group ring element $\hat{D}$ satisfies the equation
\begin{equation}\label{eq:difset}
\hat{D}\hat{D}^{(-1)}= (k-\lambda)1_{G}+\lambda \hat{G}.
\end{equation}

Two difference sets $D_{1}, D_{2} \subseteq G$ are equivalent if there is an element $g \in G$ and an automorphism $\varphi$ of $G$ such that $D_{1} = \{g \varphi(d) : d \in D_2\}$.
Difference sets are inequivalent if either they are subsets of nonisomorphic groups or if they are subsets in a common group $G$ but are not equivalent in $G$.
When we refer to finding all difference sets in a group $G$ we typically mean finding a collection of difference sets in $G$ that contains exactly one representative from each equivalence class in the complete collection of all difference sets in $G$.

If a group $G$ has a difference set $D$ then $\{gD : g \in G \}$ is the set of blocks of a symmetric $(v, k, \lambda)$ design with point set $G$.
On this design $G$ acts by left multiplication as a sharply transitive automorphism group.
Conversely, any symmetric design with a sharply transitive automorphism group on points is isomorphic to a design constructed from the set of left translates of a difference set.
While equivalent difference sets always give rise to isomorphic designs, inequivalent difference sets may also give rise to isomorphic designs, even if the difference sets belong to nonisomorphic groups.
A symmetric design is said to be genuinely non-abelian if it has no abelian group acting sharply transitively on the points of the design.
For more details on symmetric designs and difference sets, the reader may consult \cite{BeJuLe1999, IoninShrikhande2006, Jungnickel1992, Lander1983}.

Difference sets with parameters $( q^{d+1}( \frac{q^{d+1}-1}{q-1}+1), q^{d} \frac{q^{d+1}-1}{q-1}, q^{d} \frac{q^{d}-1}{q-1})$, where $q=p^m$ is a prime power, are known as McFarland difference sets.
For further discussion on McFarland difference sets, see \cite{Dillon1985, McFarland1973}.
For $q = 4$ and $d = 1$ we obtain the $(96, 20, 4)$ parameters.

There are 231 groups of order 96.
It is known that exactly 94 of these groups admit $(96, 20, 4)$ difference sets \cite{GoVuMa2005}.
Previously, all $(96, 20, 4)$ difference sets were constructed in the $75$ groups which have normal subgroups of both order $3$ and order $4$ \cite{AbuGhneimSmith2007}.
In this paper we complete the construction of all $(96, 20, 4)$ difference sets in the remaining $19$ groups.
So, all $(96, 20, 4)$ difference sets are known.
The 637 difference sets in these 19 groups provide 197 nonisomorphic symmetric $(96, 20, 4)$ designs, which brings the total number of inequivalent $(96, 20, 4)$ difference sets to 2627 and nonisomorphic symmetric $(96, 20, 4)$ designs from these difference sets to 583.
Twenty of the 583 designs can be constructed from an abelian group.
The remaining 563 designs are genuinely non-abelian.
We also use the complete collection of $(96, 20, 4)$ difference sets to construct all $(96, 20, 4, 4)$ and $(96, 19, 2, 4)$ partial difference sets and briefly look at the  corresponding strongly regular graphs with parameters $(96, 20, 4, 4)$ and $(96, 19, 2, 4)$, of which there are 58 and 12, respectively, up to isomorphism.
All of the strongly regular graphs arise from non-abelian groups.

\section{Summary of the literature on $(96,20,4)$ difference sets}\label{sec:lit}

In this paper we will refer to groups as they appear in the {\tt SmallGroups} library of the software package {\tt GAP} \cite{GAP4}.
For instance, when we work in groups of order $96$ and write {\tt [96, 14]} we mean group number 14 of order 96 in the {\tt GAP} library.

A large number of authors have contributed to the search for $(96, 20, 4)$ difference sets over the last few decades.
Here we give a summary of the results known before this paper by restating and updating the summary from \cite{AbuGhneimSmith2007}.

Of the 231 groups of order 96, only 7 are abelian.
A result of Turyn \cite{Turyn1965} rules out the existence of difference sets in the abelian groups {\tt [96, 2]} $\cong \ZZ_{96}$ and {\tt [96, 59]} $\cong \ZZ_2 \oplus \ZZ_{48}$.
Arasu, Davis, Jedwab, Ma, and McFarland \cite{ArDaJeMaMc1996} ruled out the existence of difference sets in the two additional abelian groups {\tt [96, 46]} $\cong \ZZ_4 \oplus \ZZ_{24}$ and {\tt [96, 176]} $\cong \ZZ_2 \oplus  \ZZ_2 \oplus \ZZ_{24}$.

For $q$ a prime power and $d$ a positive integer, McFarland \cite{McFarland1973} constructed difference sets with parameters $(q^{d+1}( \frac{q^{d+1}-1}{q-1}+1), q^{d}\frac{q^{d+1}-1}{q-1}, q^{d} \frac{q^{d}-1}{q-1})$ in abelian groups $G$ with an elementary abelian subgroup $H$ of order $q^{d+1}$.
McFarland's construction used hyperplanes in the projective space of dimension $d$ over $GF(q)$ to construct a special set (\emph{spread}) of $\frac{q^{d+1}-1}{q-1}$ subsets, each of size $q^d$, sitting in  $H$.
These subsets were then distributed across cosets of $H$ in $G$.
McFarland's construction with $q=4$ and  $d=1$ takes five subsets $H_1 ,H_2, \dots, H_5$ of the elementary abelian group $H \cong \ZZ_2^4$ and spreads them across five of six cosets of $H$ in $G$.
This construction gives $(96, 20, 4)$ difference sets in the two abelian groups {\tt [96, 231]} $\cong \ZZ_{2}^4 \oplus \ZZ_{6}$ and {\tt [96, 220]} $\cong \ZZ_{2}^3 \oplus \ZZ_{12}$.
Arasu and Sehgal \cite{ArasuSehgal1995} then finished the existence question for abelian groups of order 96 by constructing a $(96,20,4)$ difference set in the last abelian group, {\tt [96, 161]} $\cong \ZZ_{2} \oplus \ZZ_{4} \oplus \ZZ_{12}$.

Dillon \cite{Dillon1985} generalized McFarland's construction to work for a larger set of groups, including non-abelian groups. 
He constructed McFarland difference sets in groups that have an elementary abelian normal subgroup of order $q^{d+1}$ in their center.
Dillon's construction gives $(96, 20, 4)$ difference sets in the non-abelian groups {\tt [96, 218]} and {\tt [96, 230]}.

Undergraduate students Nichols \cite{Nichols2000} (under the supervision of Harriet Pollatsek, Mt. Holyoke College) and Axon and Gotman \cite{AxonGotman2002} (under the supervision of  Emily Moore, Grinnell College) used $(16, 6, 2)$ difference sets to construct images of $(96, 20, 4)$ difference sets in groups of order $32$ and then used those images to construct $(96, 20, 4)$ difference sets in {\tt [96, 221]} and {\tt [96, 231]}.

In his doctoral dissertation, AbuGhneim \cite{AbuGhneim2005} concentrated on generalizing previous work to non-abelian groups.
AbuGhneim and Smith \cite{AbuGhneim2004, AbuGhneim2005} constructed $(96, 20, 4)$ difference sets in groups that have $\ZZ_2^4$ as a normal subgroup.
There are $19$ such groups, which are {\tt [96, i]} where {\tt i} is in the set 
\begin{numberlists}
\numlist{70, 159, 160, 162, 167, 194, 195, 196, 197, 218, 219, 220, 221, 226, 227, 228, 229, 230, 231}.
\end{numberlists}
In addition, AbuGhneim and Smith ruled out any group $G$ which has $\ZZ_2 \oplus \ZZ_{24}$, $\ZZ_2 \oplus \mathbb{D}_{24}$, $(\ZZ_3 \rtimes \ZZ_8) \oplus \ZZ_2$, or $\mathbb{D}_{48}$ as a factor group.
This result rules out difference sets in {\tt [96, i]} where {\tt i} is in the set 
\begin{numberlists}
\numlist{6, 7, 8, 9, 11, 18, 19, 25, 28, 37, 46, 48, 55, 59, 60, 76, 80, 81, 82, 89, 93, 102, 104, 109 , 110, 111, 112, 115, 116, 127, 132, 134, 137, 176, 207}.
\end{numberlists}

In a later paper, AbuGhneim and Smith \cite{AbuGhneimSmith2007} enumerated all $(96, 20, 4)$ difference sets in groups that have normal subgroups of both order $3$ and order $4$ by using the software {\tt GAP} to build images of hypothetical $(96, 20, 4)$ difference sets in groups of order $32$ and $24$ and exhaustively searching their preimages.
There are $195$ groups of order $96$ that have normal subgroups of both order $3$ and order $4$.
Among these $195$ groups, $75$ admit $(96, 20, 4)$ difference sets and the other $120$ groups do not admit any $(96, 20, 4)$ difference sets.

Golemac, Mandi\'{c}, and Vu\v{c}i\v{c}i\'{c} \cite{GoVuMa2005, GoMaVu2007} constructed $(96,20,4)$ difference sets in $27$ non-abelian groups.
These groups are {\tt [96, i]} where {\tt i} is in the set 
\begin{numberlists}
\numlist{13, 41, 64, 70, 71, 72, 78, 87, 144, 147, 159, 160, 167, 174, 185, 186, 188, 190, 194, 195, 196, 197, 209, 226, 227, 228, 229}.
\end{numberlists}
This work was extended in \cite{GolemacVucicic}, where difference sets were constructed assuming that the group $G$ has an automorphism group isomorphic to a direct product of an abelian group of order 16 and $\ZZ_3$.
Finally, Brai\'{c}, Golemac, Mandi\'{c}, and Vu\v{c}i\v{c}i\'{c} \cite{BrGoMaVu2010} finished resolving the existence problem of $(96, 20, 4)$ difference sets in all groups of order $96$ by constructing a $(96, 20, 4)$ difference set in the last group, {\tt [96, 68]}.

Hence, a group $G$ of order $96$ admits $(96, 20, 4)$ difference sets if and only if $G$ is one of the groups {\tt [96, i]} where {\tt i} is a member of the set
\begin{numberlists}
\numlist{10, 13, 14, 20, 41, 51, 52, 54, 64, 68, 70, 71, 72, 75, 77, 78, 79, 83, 84, 85, 86, 87, 88, 90, 91, 92, 94, 95, 96, 97, 98, 99, 101, 103, 105, 129, 130, 131, 133, 135, 136, 141, 142, 143, 144, 145, 146, 147, 151, 152, 159, 160, 161, 162, 164, 165, 166, 167, 168, 169, 170, 171, 172, 173, 174, 175, 177, 185, 186, 188, 190, 191, 194, 195, 196, 197, 202, 205, 206, 209, 210, 212, 218, 219, 220, 221, 223, 225, 226, 227, 228, 229, 230, 231}.
\end{numberlists}
The number of groups of order $96$ that admit $(96, 20, 4)$ difference sets is $94$.
In \cite{AbuGhneimSmith2007} all $(96,20,4)$ difference sets up to equivalence were constructed in $75$ of these groups (the groups with normal subgroups of both order $3$ and order $4$).
In this paper, we find all $(96,20,4)$ difference sets in the remaining $19$ groups, which are {\tt [96, i]} for {\tt i} in
\begin{numberlists}
\numlist{64, 68, 70, 71, 72, 185, 186, 188, 190, 191, 194, 195, 196, 197, 202, 226, 227, 228, 229}.
\end{numberlists}
Our techniques and programs can be used to verify all previous results as well.

\section{All $(96,20,4)$ difference sets}

Let $D$ be a difference set in $G$.
A homomorphism $f$ from $G$ onto $G'$ induces, by linearity, a homomorphism from $\ZZ[G]$ onto $\ZZ[G'].$
If the kernel of $f$ is the subgroup $U$, let $T$ be a complete set of distinct representatives of cosets of $U$ and, for $g \in T$, set $t_{g} = |gU \cap D|$.
The multiset $\{t_{g} : g \in T\}$ is the collection of intersection numbers of $D$ with respect to $U$.
The image of $\hat{D}$ under the function $f$ is
\[
f(\hat{D}) = \sum_{g \in T} t_{g} f(g).
\]
By applying $f$ to both sides of Equation~\ref{eq:difset} from Section~\ref{sec:intro} we can show that $f(\hat{D})$ satisfies the equation
\begin{equation}\label{eq:contracted}
f(\hat{D})f(\hat{D})^{(-1)}= (k-\lambda)1_{G'} + \lambda |U| \hat{G'}
\end{equation}
in the group ring $\ZZ[G']$.

This contraction of $\hat{D}$ to a smaller homomorphic image often provides useful information on the existence of a difference set in the original group.
In particular, since any difference set $D$ in $G$ will have images satisfying Equation~\ref{eq:contracted} in any homomorphic image, we can first enumerate possible images of difference sets in homomorphic images of $G$ and then only search subsets of $G$ that are preimages.
The collection of homomorphic images to check using Equation~\ref{eq:contracted} and the resulting collection of possible preimages to check using Equation~\ref{eq:difset} is much smaller than the total number of subsets of $G$, which makes the search feasible.

Recall from Section~\ref{sec:lit} that $94$ groups of order $96$ admit $(96,20,4)$ difference sets and in $75$ of these groups all $(96, 20, 4)$ difference sets have been constructed by enumerating homomorphic images in groups of order 32 and 24.
The other $19$ groups that admit $(96, 20, 4)$ difference sets do not have normal subgroups of order $3$.
However, any one of these $19$ groups has either a normal subgroup of order $2$ or a normal subgroup of order $4$, which means these groups have homomorphic images of order $48$ or $24$.
To find all $(96, 20, 4)$ difference sets in these $19$ groups, we first construct all possible images in groups of order $48$ and $24$ that satisfy Equation~\ref{eq:contracted}, and then we search all preimages.

The computer algebra system {\tt GAP} was used to perform all of these computations.
This approach proves to be very efficient in finding $(96, 20, 4)$ difference sets, and enables us to construct all $(96, 20, 4)$ difference sets up to equivalence in every group of order 96.
A complete list of the 2627 inequivalent $(96, 20, 4)$ difference sets and details on which of them provide nonisomorphic symmetric $(96, 20, 4)$ designs is provided in \cite{AbuGhneimWeb}, and a summary is included in the appendix to this paper.
The {\tt DESIGN} \cite{DESIGN} package for {\tt GAP} was used to determine the 583 nonisomorphic symmetric $(96, 20, 4)$ designs that arise from these $(96, 20, 4)$ difference sets.

In addition, a general application of this homomorphic image strategy is implemented in the {\tt DifSets} \cite{DifSets} package for {\tt GAP}.
Given a group $G$, the package's main algorithm first produces a chief series
\[
G = N_1 \triangleright N2 \triangleright \dots \triangleright N_r = \{1\}
\]
of $G$ and a list of possible values of $k$ for a difference set in $G$.
For each $k$ the only possible image in $G/N_1 = \{1\}$ of any size $k$ difference set in $G$ is $k \cdot 1_{G/N_1}$.
From this starting point, the algorithm generates and checks preimages in $G/N_2$, then preimages in $G/N_3$, and so on until generating and checking preimage difference sets in $G$.
By repeatedly applying the homomorphic image idea, the search space is decreased in each step and difference sets can be efficiently enumerated for many groups, including all groups of order $96$.
The {\tt DifSets} package was used to produce the final enumerations for {\tt [96, 230]} and {\tt [96, 231]} in this paper, provides a simple interface for loading these and many other difference sets into {\tt GAP}, and can be used to recompute and verify all given results.

\section{Partial difference sets and strongly regular graphs}

Many $(96, 20, 4)$ difference sets can be used to construct partial difference sets and strongly regular graphs.
In fact, with our complete collection of all $(96, 20, 4)$ difference sets we can produce all $(96, 20, 4, 4)$ and $(96, 19, 2, 4)$ partial difference sets and generate the corresponding strongly regular graphs.
In this section we summarize the knowledge of these partial difference sets.

To start, a $(v, k, \lambda, \mu)$ partial difference set is a subset $T$ of size $k$ in a group $G$ of order $v$ such that the multiset $\{ xy^{-1} : x, y  \in T \text{ and } x \neq y \}$ contains each nonidentity element of $T$ exactly $\lambda$ times and each nonidentity element of $G \setminus T$ exactly $\mu$ times.
Two partial difference sets $T_{1}$ and $T_{2}$ in a group $G$ are equivalent if there is an automorphism $\varphi$ of $G$ such that $T_{1} = \{ \varphi(t) : t \in T_{2}\}$.

By definition, $(v, k, \lambda)$ difference sets and $(v, k, \lambda, \lambda)$ partial difference sets are identical.
However, the equivalence of partial difference sets is slightly stronger than the equivalence of difference sets.
Still, every $(96, 20, 4, 4)$ partial difference set in a group $G$ is equivalent to a translate $\{gd : d \in D\}$ for $g \in G$ of one of our collection of difference sets $D$ up to equivalence in $G$.
By enumerating all translates and removing all but one partial difference set up to equivalence we can produce the collection of all $(96, 20, 4, 4)$ partial difference sets up to equivalence.
The result is a total of 145300 $(96, 20, 4, 4)$ partial difference sets in the 94 groups of order 96 that contain difference sets.
A count of $(96, 20, 4, 4)$ partial difference sets for each group is listed in the appendix.

A subset $T$ of a group $G$ is called reversible if $T=T^{(-1)}$.
A reversible partial difference set is called regular if it does not contain the identity element.
The following results can be found in \cite{Ma1994}.

\begin{proposition}\label{prop:lambdamu}
If $T$ is a $(v, k, \lambda, \mu)$ partial difference set with $\lambda \neq \mu$, then $T$ is reversible.
\end{proposition}

\begin{proposition}
Suppose $T \subseteq G$ is a reversible set containing the identity element $1_G$.
Then $T$ is a  $(v, k, \lambda, \mu)$ partial difference set if and only if $T \setminus \{1_G\}$ is a regular $(v, k-1, \lambda-2, \mu)$ partial difference set.
\end{proposition}

\begin{proposition}
Suppose that $D$ is a $(v, k, \lambda)$ difference set in a group $G$ and $g \in G$.
Then $gD$ is a regular $(v, k, \lambda, \lambda)$ partial difference set if and only if $g^{-1} \notin D$ and $gD$ is a reversible set.
Also $gD \setminus \{1_G\}$ is a regular $(v, k-1, \lambda-2, \lambda)$ partial difference set if and only if $g^{-1} \in D$ and $gD$ is a reversible set.
\end{proposition}

From these propositions, we see that a  reversible $(96, 20, 4)$ difference set $D$ either gives a $(96, 20, 4, 4)$ regular partial difference set if $1_G \not\in D$ or a $(96, 19, 2, 4)$ regular partial difference set if $1_G \in D$.
By testing each translate of the $(96, 20, 4)$ difference sets in our collection for reversibility we can thus produce all $(96, 20, 4, 4)$ and $(96, 19, 2, 4)$ regular partial difference sets up to equivalence.
Furthermore, Proposition~\ref{prop:lambdamu} and a simple counting argument show that any $(96, 19, 2, 4)$ difference set must be regular, and so our collection of all $(96, 19, 2, 4)$ regular partial difference sets is also a collection of all $(96, 19, 2, 4)$ partial difference sets.
These regular partial difference sets are much rarer in groups of order 96 than difference sets, with a total of only 150 $(96, 20, 4, 4)$ and 33 $(96, 19, 2, 4)$ regular partial difference sets appearing in 9 groups of order 96, as listed in the appendix.

Regular partial difference sets are closely related to strongly regular graphs.
A graph is a $(v, k, \lambda, \mu)$ strongly regular graph if it has $v$ vertices where every vertex has  valency $k$, any pair of adjacent vertices have exactly $\lambda$ common neighbors, and any pair of nonadjacent vertices have exactly $\mu$ common neighbors.
We have the following known theorem that relates strongly regular graphs and partial difference sets, see \cite{BeJuLe1999}.

\begin{definition}
For a group $G$ and a subset $T$ of $G$ with $1_G \notin T$ and $T=T^{(-1)}$, the \emph{Cayley graph} $\Gamma = \mathrm{Cay}(G, T)$ is a graph whose vertex set is $G$ and two vertices $x$ and $y$ are adjacent if $xy^{-1} \in T$.
\end{definition}

\begin{theorem}
A Cayley graph $\mathrm{Cay}(G, T)$ is a $(v, k, \lambda, \mu)$ strongly regular graph if and only if $T$ is a $(v, k, \lambda, \mu)$ regular partial difference set in $G$.
\end{theorem}

Using the {\tt GRAPE} \cite{GRAPE} package for {\tt GAP} we can produce a strongly regular graph for each of the regular partial difference sets in our collection.
The final result up to isomorphism is a total of 58 $(96, 20, 4, 4)$ strongly regular graphs and 12 $(96, 19, 2, 4)$ strongly regular graphs.
A summary is included in the appendix, with more details at \cite{AbuGhneimWeb}.

The study of $(96, 20, 4, 4)$ and $(96, 19, 2, 4)$ partial difference sets has a long history.
In 1969, Ahrens and Szekeres \cite{AhrensSzekeres1969} constructed a generalized quadrangle, $GQ(5, 3)$, with point graph $(96, 20, 4)$.
The line graph of the generalized quadrangle $GQ(5, 3)$ is a strongly regular graph with parameters $(96, 20, 4, 4)$.
That graph has an automorphism group of order 138240 and includes several subgroups of order 96 acting sharply transitively on points. 
However, at the time this result was not viewed in terms of strongly regular Cayley graphs or partial difference sets.

A summary of the Eighth International Conference on Geometry, University of Haifa, March 7-14, 1999 references a talk by Klin, ``Strongly regular Cayley graphs on 96 vertices".
That article is the first report (that we know of) on a Cayley graph or partial difference set with parameters $(96, 20, 4, 4)$.
The graph was constructed using the software package {\tt COCO} \cite{FaradzevKlin1991}.

In \cite{BrKoKl2003}, Brouwer, Koolen, and Klin reported on another $(96, 20, 4, 4)$ strongly regular graph constructed by examining a rooted graph on 192 vertices, using it to construct a distance regular graph of diameter three on 96 vertices and then merging the classes of distance 1 and 3 to create the strongly regular graph.
This graph also has a large automorphism group (of order 11520) and several subgroups of order 96 acting sharply transitively on the points.
The authors then use the same root graph to find a second nonisomorphic (96, 20, 4, 4) strongly regular graph with the same automorphism group.
This paper says, ``No doubt there are lots of graphs and designs with these parameters."
{\em Yes!}
In \cite{KlinReichard2003}, Klin and Reichard follow up on \cite{BrKoKl2003} by using the concept of partial linear space to generalize generalized quadrangles and explain the two earlier strongly regular graphs in terms of partial linear spaces.

More recently, Golemac, Mandi\'{c} and Vu\v{c}i\v{c}i\'{c} \cite{GoVuMa2005} found that $9$ groups of order $96$ admit regular $(96, 19, 2, 4)$ and $(96, 20, 4, 4)$ partial difference sets.
These groups are {\tt [96, i]} where {\tt i} is in
\begin{numberlists}
\numlist{64, 70, 71, 186, 190, 195, 197, 226, 227}.
\end{numberlists}
These are exactly the 9 groups in which such partial difference sets exist, as verified by our exhaustive search.
They found $29$ inequivalent regular $(96, 19, 2, 4)$ and $115$ inequivalent regular $(96, 20, 4, 4)$ partial difference sets, and provide more details on these partial difference sets on the webpage \cite{GoMaVuWeb}.

Furthermore, in \cite{GoMaVu2006} Golemac, Mandi\'{c}, and Vu\v{c}i\v{c}i\'{c} use lists of symmetric designs generated earlier along with {\tt GAP} and {\tt GRAPE} to search for reversible difference sets, finding them in the same nine groups as above.
This creates six $(96, 20, 4, 4)$ graphs and two $(96, 19, 2, 4)$ graphs.
Two of the graphs on 96 vertices have full automorphism group of order 96 (specifically, {\tt [96, 195]}).
The graph with largest automorphism group (of order 138240) is the collinearity group of the generalized quadrangle $GQ(5, 3)$.
The graph they call $\Gamma_2$ with the next largest automorphism group (of order 11520), is the one found by Brouwer, Koolen, and Klin \cite{BrKoKl2003}.

Law, Praeger, and Reichard \cite{LaPrRe2009} give four symmetric $2-(96,20,4)$-designs with flag transitive automorphism groups, the three occurring in \cite{BrKoKl2003} and one more.
Each of these gives a strongly regular graph.
There are four flag transitive symmetric designs \cite{Reichard2003, LaPrRe2009}.

Our complete enumeration of all $(96, 20, 4)$ difference sets thus provides an additional $4$ inequivalent regular $(96, 19, 2, 4)$ and another $35$ inequivalent regular $(96, 20, 4, 4)$ partial difference sets, finally giving a definitive answer to the collections of these structures and their corresponding strongly regular graphs.

\section*{Appendix}

The following tables list the counts of difference sets and related structures found in each group of order 96.
Each line is a group listed by its index {\tt [96, i]} in the {\tt SmallGroups} library.
Groups that do not appear in the tables contain no difference sets or related structures.
A list of the actual difference sets and details on the structures they generate can be found in \cite{AbuGhneimWeb}, and some tools to recompute these results can be found in \cite{DifSets}.

The first table contains counts of difference sets, designs, and partial difference sets for the 94 groups that contain difference sets.
All counts are up to equivalence or isomorphism, but note that some symmetric designs from nonisomorphic groups are isomorphic.
The total number of symmetric $(96, 20, 4)$ designs up to isomorphism found from these groups is 583.

The second table contains counts of regular partial difference sets and strongly regular graphs for the 9 groups that contain regular partial difference sets.
Note that all $(96, 19, 2, 4)$ partial difference sets are regular, so this is a complete count of the $(96, 19, 2, 4)$ partial difference sets.
As with designs, some strongly regular graphs from nonisomorphic groups are isomorphic.
The total number of $(96, 20, 4, 4)$ strongly regular graphs and $(96, 19, 2, 4)$ strongly regular graphs up to isomorphism found using these groups is 58 and 12, respectively.

\begin{center}
\begin{longtable}{cccc}
\multicolumn{1}{p{1.5cm}}{\centering \hfill\break \hfill\break Group} &
\multicolumn{1}{p{3cm}}{\centering \hfill\break $(96, 20, 4)$ \\ Difference \\ Sets} &
\multicolumn{1}{p{3cm}}{\centering \hfill\break $(96, 20, 4)$ \\ Symmetric \\ Designs} &
\multicolumn{1}{p{3cm}}{\centering $(96, 20, 4, 4)$ \\ Partial \\ Difference \\ Sets} \\
\hline
{\tt [96, 10]} & 4 & 3 & 216 \\
{\tt [96, 13]} & 16 & 15 & 832 \\
{\tt [96, 14]} & 4 & 4 & 384 \\
{\tt [96, 20]} & 8 & 8 & 768 \\
{\tt [96, 41]} & 16 & 15 & 832 \\
{\tt [96, 51]} & 8 & 4 & 254 \\
{\tt [96, 52]} & 4 & 4 & 384 \\
{\tt [96, 54]} & 12 & 8 & 816 \\
{\tt [96, 64]} & 14 & 10 & 620 \\
{\tt [96, 68]} & 2 & 2 & 132 \\
{\tt [96, 70]} & 28 & 22 & 1012 \\
{\tt [96, 71]} & 8 & 4 & 416 \\
{\tt [96, 72]} & 2 & 2 & 132 \\
{\tt [96, 75]} & 88 & 84 & 4576 \\
{\tt [96, 77]} & 40 & 32 & 2080 \\
{\tt [96, 78]} & 10 & 10 & 488 \\
{\tt [96, 79]} & 24 & 24 & 1248 \\
{\tt [96, 83]} & 18 & 14 & 904 \\
{\tt [96, 84]} & 72 & 72 & 3744 \\
{\tt [96, 85]} & 136 & 96 & 7072 \\
{\tt [96, 86]} & 48 & 40 & 2496 \\
{\tt [96, 87]} & 16 & 10 & 832 \\
{\tt [96, 88]} & 48 & 36 & 2496 \\
{\tt [96, 90]} & 24 & 16 & 1248 \\
{\tt [96, 91]} & 8 & 4 & 416 \\
{\tt [96, 92]} & 40 & 24 & 2080 \\
{\tt [96, 94]} & 84 & 64 & 4368 \\
{\tt [96, 95]} & 96 & 72 & 4992 \\
{\tt [96, 96]} & 120 & 60 & 6240 \\
{\tt [96, 97]} & 16 & 16 & 832 \\
{\tt [96, 98]} & 24 & 20 & 1248 \\
{\tt [96, 99]} & 12 & 12 & 624 \\
{\tt [96, 101]} & 8 & 4 & 416 \\
{\tt [96, 103]} & 40 & 32 & 2080 \\
{\tt [96, 105]} & 8 & 8 & 416 \\
{\tt [96, 129]} & 36 & 36 & 2224 \\
{\tt [96, 130]} & 88 & 74 & 5984 \\
{\tt [96, 131]} & 40 & 36 & 2080 \\
{\tt [96, 133]} & 16 & 16 & 832 \\
{\tt [96, 135]} & 12 & 8 & 624 \\
{\tt [96, 136]} & 8 & 8 & 416 \\
{\tt [96, 141]} & 60 & 44 & 3120 \\
{\tt [96, 142]} & 64 & 52 & 3328 \\
{\tt [96, 143]} & 12 & 12 & 624 \\
{\tt [96, 144]} & 8 & 6 & 416 \\
{\tt [96, 145]} & 4 & 4 & 208 \\
{\tt [96, 146]} & 32 & 20 & 1664 \\
{\tt [96, 147]} & 4 & 2 & 208 \\
{\tt [96, 151]} & 36 & 36 & 1872 \\
{\tt [96, 152]} & 12 & 8 & 624 \\
{\tt [96, 159]} & 20 & 18 & 1392 \\
{\tt [96, 160]} & 30 & 25 & 1580 \\
{\tt [96, 161]} & 6 & 6 & 280 \\
{\tt [96, 162]} & 28 & 28 & 2512 \\
{\tt [96, 164]} & 20 & 20 & 1040 \\
{\tt [96, 165]} & 16 & 16 & 832 \\
{\tt [96, 166]} & 16 & 16 & 832 \\
{\tt [96, 167]} & 56 & 55 & 4140 \\
{\tt [96, 168]} & 24 & 24 & 1600 \\
{\tt [96, 169]} & 8 & 8 & 416 \\
{\tt [96, 170]} & 44 & 44 & 2288 \\
{\tt [96, 171]} & 16 & 16 & 1008 \\
{\tt [96, 172]} & 16 & 16 & 832 \\
{\tt [96, 173]} & 16 & 16 & 832 \\
{\tt [96, 174]} & 6 & 6 & 456 \\
{\tt [96, 175]} & 8 & 8 & 416 \\
{\tt [96, 177]} & 16 & 16 & 1536 \\
{\tt [96, 185]} & 20 & 15 & 896 \\
{\tt [96, 186]} & 16 & 11 & 512 \\
{\tt [96, 188]} & 52 & 29 & 3712 \\
{\tt [96, 190]} & 40 & 24 & 2560 \\
{\tt [96, 191]} & 8 & 8 & 768 \\
{\tt [96, 194]} & 72 & 31 & 3264 \\
{\tt [96, 195]} & 84 & 29 & 4720 \\
{\tt [96, 196]} & 82 & 53 & 4456 \\
{\tt [96, 197]} & 72 & 39 & 4288 \\
{\tt [96, 202]} & 25 & 25 & 2400 \\
{\tt [96, 205]} & 48 & 48 & 2848 \\
{\tt [96, 206]} & 8 & 8 & 416 \\
{\tt [96, 209]} & 4 & 2 & 208 \\
{\tt [96, 210]} & 20 & 20 & 1040 \\
{\tt [96, 212]} & 16 & 16 & 832 \\
{\tt [96, 218]} & 14 & 14 & 760 \\
{\tt [96, 219]} & 4 & 4 & 272 \\
{\tt [96, 220]} & 12 & 12 & 976 \\
{\tt [96, 221]} & 12 & 12 & 976 \\
{\tt [96, 223]} & 12 & 12 & 624 \\
{\tt [96, 225]} & 6 & 6 & 314 \\
{\tt [96, 226]} & 28 & 17 & 836 \\
{\tt [96, 227]} & 42 & 30 & 1672 \\
{\tt [96, 228]} & 34 & 32 & 1528 \\
{\tt [96, 229]} & 8 & 8 & 408 \\
{\tt [96, 230]} & 2 & 2 & 32 \\
{\tt [96, 231]} & 2 & 2 & 72 \\
\end{longtable}
\end{center}

\begin{center}
\begin{longtable}{ccccc}
\multicolumn{1}{p{1.2cm}}{\centering \hfill\break \hfill\break Group} &
\multicolumn{1}{p{2.5cm}}{\centering $(96, 20, 4, 4)$ \\ Regular \\ Partial \\ Difference \\ Sets} &
\multicolumn{1}{p{2.5cm}}{\centering $(96, 19, 2, 4)$ \\ Regular \\ Partial \\ Difference \\ Sets} &
\multicolumn{1}{p{2.5cm}}{\centering $(96, 20, 4, 4)$ \\ Strongly \\ Regular \\ Graphs} &
\multicolumn{1}{p{2.5cm}}{\centering $(96, 19, 2, 4)$ \\ Strongly \\ Regular \\ Graphs} \\
\hline
{\tt [96, 64]} & 7 & 2 & 6 & 2 \\
{\tt [96, 70]} & 10 & 2 & 10 & 2 \\
{\tt [96, 71]} & 7 & 1 & 3 & 1 \\
{\tt [96, 186]} & 14 & 2 & 13 & 2 \\
{\tt [96, 190]} & 8 & 2 & 8 & 2 \\
{\tt [96, 195]} & 48 & 12 & 41 & 8 \\
{\tt [96, 197]} & 14 & 2 & 13 & 2 \\
{\tt [96, 226]} & 23 & 5 & 22 & 5 \\
{\tt [96, 227]} & 19 & 5 & 19 & 4 \\
\end{longtable}
\end{center}

\bibliographystyle{plain}
\bibliography{AbPeSm2019}

\end{document}